# The principle of maximum in the imitative control tasks

Ivan V. Kazachkov[1,2]

Dept of Information Technologies and Data Analysis, Nizhyn Gogol state university, Ukraine
Dept of Energy Technology, Royal Institute of Technology (KTH), Stockholm, Sweden)
Ivan.Kazachkov@energy.kth.se

**Abstract:** The article is devoted to the problem of applying the maximum principle for finding optimal control parameters in simulation tasks of interest for a variety of engineering and industrial systems and processes. Especially important is the problem for such systems where it is practically impossible to organize a control system because of an unknown model or because it is impossible to strictly follow the specified trajectories of control parameters in each particular case. A so-called $\Delta$-procedure is constructed that converges and allows one to obtain an approximation and optimal control of the system in a finite number of steps with a given accuracy $\Delta$. The theorem on the $\Delta$-procedure is proved and illustrative examples are given. Some features of the application of the proposed algorithm to real problems of simulation control are discussed, e.g. the optimal number of points dividing the time interval under consideration, and tasks for further individual research.

## 1. Introduction to the problem

Many objects of a control are so complex that it is impossible to create the precise advanced models for their use in the control systems. Therefore, development of the adaptive control systems, which could imitate the behaviors of some etalon system, after their education, is important problem in a modern control theory [1-5]. In a number of the industrial, technological and other control systems, the ergatic control systems are applied, which are operating together with a team of the experienced operators supplying satisfactory running of the complex systems. For example, the steel making processes are operated with an application of the imitative control using the data base of the successful "trajectories" of the controlling process. The system of adaptive control is built on the principle of the best closeness of the imitative trajectory to the one from the data base of the optimal trajectories at each reference moment of time. The sophisticated machine control has a learning unit, which is analyzing the data base and choosing the best trajectories from it based on the stated criterion. Thus, from this point of view, it is the imitative control system. Presently the machine learning is a growing direction of the theory and applications. It includes slightly wider field as far as it is applied to the many other problem solving tasks including some kind of a study helper.

## 2. Statement of the task by imitative control

In the task of the imitative control, in a creation of such systems, the peculiarity consists in an absence of the mathematical model of the control object. Thus, the mathematical model (e.g. equations of the object) and elements of the control system are unknown. This substantially complicates a solution of the problem. Therefore, a compromise is foreseen both, in the strict statement of the problem, as well as in the accuracy of its solution. The problem on optimal control of such objects is especially important because the optimal trajectories may be absent among the positive ones in the data base. Unfortunately, the task of optimal control in the complex nonlinear objects sometimes is not resolved even in a case when the mathematical model for it is known.

Accounting the above-mentioned, we consider in continuation of the earlier work [1] the following task on optimal control of the object of unknown structure, which is known by its positive and negative operating information. We consider the information about the controlling object as positive if the object works in its normal regime, otherwise - negative if the working regime does not satisfy the requirements. In general, the mathematical model of the object under control can be presented as

$$\frac{dx}{dt} = f(x,u,t), \qquad (1)$$

where $n$-dimensional vector of the parameters of state for the object, $u(t)$ is the $r$-dimensional vector of the control parameters, $t$ - time, $f$ – $n$-dimensional vector-function. The control task is stated as follows. An object (1) must be guided in a way to get from the initial state $x(t_0)$ to the given final state $x(t_1)$ for minimal time (the task on optimal speed of the process), where $t_0$, $t_1$ are the initial and final moments of the control process.

### 3. The proposed algorithm based on the principle of maximum

The following algorithm for finding the optimal control of the object (1) was proposed based on the Pontryagin's maximum principle. Having in hands the experimental data about the states of object with time, we can do as follows. The time interval is divided in $N_m$ parts (for $m$-th case): $\tau_0$, $\tau_1$, $\tau_2$,..., $\tau_{N_m}$ ($\tau_0 = t_0$, $\tau_{N_m} = t_1$), where $N_m$ is a number of the intervals for the case of $m$-th division of the interval $t \in [t_0, t_1]$ on $N_m$ subintervals. And for each of such cases the piecewise-continuous approximations of the system (1) by linear systems of the form:

$$\frac{dx_m}{dt} = A_m x_m + B_m u_m, \qquad (2)$$

are done, where $A_m(t)$, $B_m(t)$ are the piecewise-constant matrices; $x_m(t)$, $u_m(t)$ are the piecewise-continuous vectors. The matrices $A_m(t)$, $B_m(t)$ are built based on the given information (both, positive and negative) about the states of the object. The system (2) can be presented in a form:

$$\frac{dx_{mk}}{dt} = A_{mk} x_{mk} + B_{mk} u_{mk}, \qquad (3)$$

where for $\forall k = \overline{1, N_m}$ there is $A_{mk}, B_{mk} - const$,

$$x_{mk} = x_m - x_m(\tau_k). \qquad (4)$$

For the linear systems (3), (4), the optimal control is determined as follows, according to the Pontryagin's maximum principle [3]: the control parameter $u_{mk}(t)$ ($t \in [\tau_{k-1}, \tau_k]$) is managing transfer of the system from the state $x_{mk}(\tau_{k-1})$ into the state 0, being the optimal one in the sense of speed, if it is satisfying the Pontryagin's maximum principle:

$$\psi_{mk}(t) B_{mk}(t) u_{mk}(t) = \max_{u_{mk} \in U} \psi_{mk}(t) B_{mk} u_{mk}, \qquad (5)$$

where $U$ is a region of the attainable control parameters, $\psi_{mk}(t)$ is nontrivial solution of the system

$$\frac{d\psi_{mk}}{dt} = -A_{mk}^{*}\psi_{mk}. \tag{6}$$

Here $A_{mk}^{*}$ is a matrix transposed to $A_{mk}$. It is known [3] that for the systems (3), (4) the condition $\psi_{mk}(t)x_{mk}(t) = const$ is satisfied. And equation (5) determines the piecewise-continuous optimal control $u_{mk}(t)$, where according to the Feldbaum A.A. [6], each function $u_{mk}(t)$ has only extremal values and no more than $n-1$ switches ($n$ is an order of the system (1)). Increasing a quantity $N_m$ of the divisions of interval $t \in [t_0, t_1]$ we can obtain more and more approximations close to the system (1). But a number of switches will grow being limited from the top by $(n-1)N_m$. Therefore, the conditions for optimal division of the interval must be got, thus, obtained the minimal value $N_m$, which supplies the given accuracy of the approximation (1). The obtained piecewise-continuous control is optimal.

According to the Pontryagin's principle of maximum the Hamilton function $H = \psi_i f_i$, $i = \overline{1,n}$, for the system (1) has its maximal values on the optimal trajectories, therefore the condition

$$\max_{N_m} \frac{1}{t_1 - t_0} \sum_{k=1}^{N_m} H_{mk}\tau_{mk}, \tag{7}$$

must satisfy, where $\tau_{mk} = (t_1 - t_0)\varepsilon_{mk}/N_m$; $\varepsilon_{mk}$ - some dimensionless value satisfying the condition $\sum_{k=1}^{N_m}\varepsilon_{mk} = N_m$, in particular can be $\varepsilon_{mk}=1$ put. The second condition for the function $H_{mk}$ may be stated minimal mean deviation of $H_{mk}$ from the average value of this function $H_m$ on the interval $t \in [t_0, t_1]$:

$$\min_{N_m} \frac{1}{N_m} \sum_{k=1}^{N_m} (H_{mk} - H_m)\varepsilon_{mk}, \tag{8}$$

where $H_m = (1/N_m)\sum_{k=1}^{N_m} H_{mk}\varepsilon_{mk}$.

Thus, for determination of the optimal control parameters of the system (1) with an unknown right side, which transfer the system from the initial state $x(t_0)$ to the final state $x(t_1)$ for minimal time, it is required after division of the interval $[t_0, t_1]$ into a few parts, approximate the system (1) by piecewise-continuous system (2) using the experimental data and determining the piecewise-continuous control parameters. Then, increasing a number of the dividing points of the temporal interval, continuing the described process until the condition,

$$\min_{N_m}(t_1 - t_0) - \min_{N_{m-1}}(t_1 - t_0) \leq \Delta \tag{9}$$

is satisfied, where $\min(t_1-t_0)$ for the $m$-th and $(m-1)$-th divisions of the interval $[t_0, t_1]$ are computed according to the above-described. Here $\Delta$ is given attainable inaccuracy.

## 4. The proposed Δ-procedure based on the principle of maximum

The procedure (9) thus proposed was named the Δ-procedure, which was introduced by the following theorem [1]:

***Theorem***. Δ-procedure converges and allows obtaining the approximation and optimal control for the system (1) in a finite number of the steps with a given accuracy.

*Proof of theorem*. Let $N_m \to \infty$, then, as far as the left and right sides of the equations (1) and (3) coincide in all interval $t \in [t_0, t_1]$, the system (3), (4) is transferring into (1). What is more, the Hamilton functions and the optimal control functions coincide too. Supposed $\tau = (t - t_0)/(t_1 - t_0)$ in (7), we can get the condition $\max \max_{u \in U} \int H(\tau) d\tau$, which corresponds to $\max \max_{u \in U} H$. The condition (8) is satisfied automatically. Thus, Δ-procedure converges. Consecutively, for the finite number of the steps determined by required accuracy Δ, one can get the approximation and optimal control for the system (1).

In case of the general problem of optimal control for the object (1) under condition of a minimum of the quality functional

$$I = \int_{t_0}^{t_1} f_0(x, u) dt \quad (10)$$

where $f_0(x,u)$ is a function chosen by the statement, the object (1) must be transferred from the initial state $x(t_0)$ to the final state $x(t_1)$ in such way that the functional (10) gets its extreme value.

The control $u(t)$, a solution of the stated task, is optimal. If $f_0 \equiv 1$, $I = t_1 - t_0$, the task reduces to the optimal performance problem. Adding to the space of the states of object the new variable $x_0$, which follows according to (10) $dx_0/dt = f_0(x, u)$. Now introducing for each trajectory the new time $\tau$ connected to the old time by correlation $d\tau = f_0(x, u)dt$, we get $I = \tau_1 - \tau_0$, so that in a new time the task stated is transformed into the optimal performance problem. Thus, the Δ-procedure is similarly applied in a general case as well. If the problem with varying boundaries is considered, so that $x(t_0) \in M_0$, $x(t_1) \in M_1$, where the mathematical sets $M_0$, $M_1$ are the varieties (the lines or surfaces of the stated number of the measurements), then all the above-mentioned is kept. But the transversality conditions of the vector $\psi(t)$ are added according to the Pontryagin's maximum principle [3]: vector $\psi(t_0)$ is orthogonal to all tangential vectors of the manifold $M_1$ at the point $x(t_1)$.

## 5. Examples illustrating the proposed Δ-procedure

Example 1. We consider the following model of the system

$$dx/dt = x^2 + u^2, \quad (11)$$

where $x(t_0) = 0$, $x(t_1) = 1$ and $u \in [-1, 1]$. Investigating (11) by the maximum principle, we get the Hamilton function $H = \psi(x^2 + u^2)$. Here $\psi$ is computed as solution of the equation $d\psi/dt = -2x\psi$, and afterward $\max_{u \in [-1,1]} H$ is found, which yields that $u = 1$ is the optimal control for the system (11). The corresponding optimal trajectory is $x = \text{tg}\, t$, which results in a solution $\min (t_1 - t_0) = \pi/4 \approx 0.78$, which is valid only until the breaking point $t = \pi/2 \approx 1.57$. Then applying the Δ-procedure we perform a few

dividing of the interval and determine the corresponding approximations of (11) and the optimal trajectories:

$$t = \tau_0, \ x=0; \quad t = \tau_1, \ x=0.5; \quad t = \tau_2, \ x=1; \quad u=0.5. \qquad (12)$$

System (11) is approximated by piecewise-continuous systems using the known points (12):

$$t \in [\tau_0, \tau_1], \quad dx/dt = a_1 x + b_1 u, \qquad (13)$$
$$t \in [\tau_1, \tau_2], \quad dx/dt = a_2 x + b_2 u,$$

where the constants $a_j$, $b_j$ ($j=1,2$) are computed substituting the (13) into the conditions (12), with account of the corresponding values $x$ obtained from (11), (12). We have got as follows:

$$t \in [\tau_0, \tau_1], \quad a_1 = b_1 = 0.5; \qquad (14)$$
$$t \in [\tau_1, \tau_2], \quad a_2 = b_2 = -0.5.$$

The obtained approximations (13), (14) investigated according to the principle of maximum result in the following

$$t \in [\tau_0, \tau_1], \quad H_1 = 0.5(x+u)\psi_1, \quad d\psi_1/dt = -0.5 x \psi_1,$$

where $H_1 = 0.5(x+u)\exp(-0.5 \int x\,dx)$. Studying the $H_1$ on maximum by the control parameter $u \in [-1,1]$ we get the optimal control $u=1$. Similarly, the case $t \in [\tau_1, \tau_2]$ was studied. The corresponding optimal trajectories were got as follows:

$$t \in [\tau_0, \tau_1], \quad t - \tau_0 = 2\ln(x+1);$$
$$t \in [\tau_1, \tau_2], \quad t - \tau_1 = (2/3)\ln[(3x+1)/2.5].$$

Here $u=1$ according to the principle of maximum at the interval $t \in [\tau_0, \tau_1]$, and $u=-1$ at the interval $t \in [\tau_1, \tau_2]$, though an approximation was made by the points with the control value $u=0.5$. The optimal trajectories result: $\min(t_1 - t_0) \approx 0.82 + 0.31 \approx 1.13$.

Example 2. Similarly, the next cases were considered ($u=1$ is optimal control everywhere):

1) $t = \tau_0, x=0; \ t = \tau_1, x=0.5; \ t = \tau_2, x=1; \ u=0.9;$

    $t \in [\tau_0, \tau_1]$, $t - \tau_0 = 2\ln(0.5x/0.9 + 1)$;

    $t \in [\tau_1, \tau_2]$, $t - \tau_1 = (2/3)\ln[(1.5x + 0.34)/1.09]$, $\min(t_1 - t_0) \approx 1.1$.

2) $t = \tau_0, x=0; \ t = \tau_1, x=0.25; \ t = \tau_2, x=0.5; \ t = \tau_3, x=0.75; \ t = \tau_4, x=1; \ u=0.9;$

    $t \in [\tau_0, \tau_1]$, $t - \tau_0 = 4\ln(x/3.6 + 1)$;

    $t \in [\tau_1, \tau_2]$, $t - \tau_1 = (4/3)\ln[(x+1.1)/1.26]$;

    $t \in [\tau_2, \tau_3]$, $t - \tau_2 = 0.8\ln[(1.25x + 0.483)/1.008]$;

    $t \in [\tau_3, \tau_4]$, $t - \tau_3 = (20/17)\ln[(0.85x + 1.07)/1.71]$; $\min(t_1 - t_0) \approx 0.94$.

3) $t = \tau_0$, $x=0$; $t = \tau_1$, $x=0.75$; $t = \tau_2$, $x=1$; $u=1$;

$$t - \tau_0 = (4/3)\ln(0.75x+1),$$

$$t - \tau_1 = (4/7)\ln[(7x+1)/6.25], \quad \min(t_1 - t_0) \approx 0.74.$$

4) $t = \tau_0$, $x=0$; $t = \tau_1$, $x=0.5$; $t = \tau_2$, $x=1$; $u=1$;

$$t - \tau_0 = 2\ln(0.5x+1),$$

$$t - \tau_1 = (2/3)\ln[(1.5x+0.5)/1.25], \quad \min(t_1 - t_0) \approx 0.76.$$

## 6. Discussions

For the nonlinear equation (11) an optimal control was got $u=1$ and the corresponding optimal solution $x=\text{tg}\,t$, which is stable. But for some real systems described by the simplest development equation with deviating arguments

$$dx/dt = kx(t \pm \tau), \tag{14}$$

there are available unstable regimes. Here $\tau > 0$ is a time deviating term, $k$ is the coefficient (intensity of a development), positive for the system's grow and negative for decreasing of a system. The sign minus in the equation (14) corresponds to a time delay, plus – to a time forecasting. As shown in [12], in many different systems of the very wide real systems (technical and organisms), the development is going with account of both, delay and forecasting. And there are available some stable processes of a system's development only in a range from $\tau = 1/u$ till $\tau = 1.293/u$, so that solution of (14) is built as $x=x_0\exp(ut)$ with a piecewise-continuous value of the control parameter $u$.

As shown by the results above, the closer are experimental data to the optimal ones, the faster $\Delta$-procedure converges. Therefore, the positive information is better for approximation of the object (1) because it reduces the calculations substantially.

The results obtained can be used for optimization of the imitative control systems, for example the ones by a control of the steel making processes [2]. The open-hearth and electric arc processes are distinguished by the presence of measurements of output parameters at any time during the melting process (sample selection), on the basis of which the steelmaker performs correction of fusible control. Stating some value of the required inaccuracy $\Delta$, one can apply the above-described $\Delta$-procedure to obtain, based on the positive information about the object, the piecewise-continuous linear approximations for the object and the corresponding optimal control parameters. Then the algorithm of the open-hearth fusible process can be reduced to the switching of the control parameters in a time moments stated, which must be performed by the steelmaker. In this process, the main task the steelmaker is exact fulfillment of the instruction by control, which is prepared with account of the optimal control parameters and the known optimal trajectories.

An application of the $\Delta$-procedure for a control of the converter steelmaking process is complicated by the fact that the vector of trajectories is unknown [2]. But some information about the object may be presented by a vector of the indirect parameters, e.g. content of the off-gas composition, which can be used for optimization of the process by the above-described algorithm.

The problem of a system's identification, having important practical value as far as a general algorithm cannot work without a simple linear model of the object, was not considered in this paper. We have also to underline an insurmountable obstacle in solving the problem. By $m \to \infty$, the control task is solved successfully but the identification task is hardly solved, while by small values of $m$ the situation changes to the opposite one. Thus, the task on finding the optimal value of *m* is stated for a separate investigation.

These and similar algorithms can be applied in different systems of diverse nature [7-11]. Except the considered above, there are important problems on stability and critical levels of the automatically controlled systems too, especially with account of the deviated arguments, which are nearly always present in a system [10-12]. The optimal control can be unstable, especially in case of deviating arguments; therefore, determination of the optimal control must be done together with a stability analysis [11, 12]. Different approaches for an optimal control developed in [3] (Pontryagin's maximum principle) and [12] (critical levels and stable development in systems with delays and time forecasting), both reveal the piecewise-continuous control functions for the optimal development.

## 7. The conclusion

The proposed algorithm and the $\Delta$-procedure for solution of the tasks by imitative control based on application the Pontryagin's maximum principle were proven by the theorem and illustrated by a few examples. Finding of the optimal control parameters in the tasks of imitative control is of interest for many applications. For example, it is important for imitative control of the steel casting. The specific feature of an application of the proposed algorithm to the real problems of imitative control is an optimal number of the points dividing the considered time interval. This is stated as a problem for further separate studies, as well as the problem of mathematical modeling of the object and its identification.

**References**

[1] Vasil'ev V.I., Kazachkov I.V. Maximum principle in simulation control// Kibernetika i Vychislitel'naya Tekhnika.- 1983.- P. 23-28.
[2] Vasil'ev V.I., Koval' P.N., Konovalenko V.V., Lazarenko I.Yu. Principles of imitative control// Soviet automatic control.- 1982, №1, P. 60-69; №2, P. 58-65.
[3] Pontryagin L.S., Boltyanskii V.G., Gamkrelidze R.V., Mizhenko E.F. *Mathematical theory of optimal processes*. (John Wiley & sons, New-York, 1962).
[4] Pesch H.J. and Plail M. The Maximum Principle of optimal control: A history of ingenious ideas and missed opportunities// Control and Cybernetics, vol. 38 (2009), No. 4A, P. 973-995.
[5] Kim A.V., Ivanov A.V. Systems with Delays: Analysis, Control, and Computations.- Wiley, Scriver Publ., 2015, 167 pp.
[6] Feldbaum A. A. Optimal systems. Disciplines and Techniques of System Control, Chap. VII, J. Peschon, ed. Random House (Blaisdell), New York, 1965.
[7] Cross K., Torrisi S., Losin E.AR., Iacoboni M. Controlling automatic imitative tendencies: Interactions between mirror neuron and cognitive control systems// NeuroImage.- June 2013, 83.
[8] Klapper A., Ramsey R., Wigboldus D., Cross E.S. The Control of Automatic Imitation Based on Bottom-Up and Top-Down Cues to Animacy: Insights from Brain and Behavior// Journal of Cognitive Neuroscience.- April 2014, P.1-11.
[9] Li X.-Q., Zeng S.-S. Model-establishing and experimental imitation for leaning electric stove's flow and automatic control system.- June 2006 Zhongnan Daxue Xuebao (Ziran Kexue Ban)/Journal of Central South University (Science and Technology) 37(3): 547-552.
[10] Boccia A., Vinter R.B. The Maximum Principle for Optimal Control Problems with Time Delays// SIAM Journal on Control and Optimization.- January 2017.- 55(5): 2905-2935.
[11] Kazachkov I.V. Determination and Calculation the Critical Levels in Development of Complex Systems of Different Nature with Shifted Arguments for their Investigation and Optimal Control// WSEAS transactions on systems and control.- Volume 12, 2017, P. 21-30.
[12] Zhirmunsky, A.V. and Kuzmin, V.I., *Critical Levels in the Development of Natural Systems*, Berlin: Springer, 1988.